\documentclass{article}
\usepackage{amsmath}
\title{A NEW TYPE OF CONTINUED FRACTION EXPANSION}
\author{Ion COLTESCU, Dan LASCU\\
\small ``Mircea cel Batran'' Naval Academy, 1 Fulgerului,\\
\small 900218 Constanta, Romania\\
\small {\it e-mail}: icoltescu@yahoo.com, lascudan@gmail.com 
}
\date{}
\begin{document}
\maketitle
\thispagestyle{empty}
\begin{abstract}
In this paper we define a new type of continued fraction expansion for a real number $x \in I_m:=[0,m-1], m\in N_+, m\geq 2$:
\[
    x = \frac{m^{-b_1(x)}}{\displaystyle 1+\frac{m^{-b_2(x)}}{1+\ddots}}:=[b_1(x), b_2(x), \ldots]_m.
\]
Then, we derive the basic properties of this continued fraction expansion, following the same steps as in the case of the regular continued fraction expansion. The main purpose of the paper is to prove the convergence of this type of expansion, i.e. we must show that
\[
    x= \lim_{n\rightarrow\infty}[b_1(x), b_2(x), \ldots, b_n(x)]_m.
\]

{\bf Keywords:} {\it continued fractions, incomplete quotients}
\end{abstract}

\section{INTRODUCTION}
In this section we make a brief presentation of the theory of
regular continued fraction expansions.

It is well-known that the regular continued fraction expansion of a real number looks as follows:
\begin{equation}
    \frac{1}{\displaystyle a_1+\frac{1}{\displaystyle a_2+\ddots + \frac{1}{a_n + \ddots}}}
\end{equation}
where $a_n \in N$, $\forall n \in N_+$. We can write this
expression more compactly as
\[
    [0;a_1, a_2, \ldots, a_n, \ldots].
    \label{eq11}
\]
The terms $a_1, a_2, \ldots$  are called the incomplete quotients of the continued fraction.
Continued fractions theory starts with the procedure known as Euclid's algorithm for finding the greatest common divisor.
To generalize Euclid's algorithm to irrational numbers from the unit interval $I$, consider the continued fraction transformation
$\tau:I \rightarrow I$ defined by
\begin{equation}
    \tau(x):=\frac{1}{x} - \left[\frac{1}{x}\right], x \neq 0, \tau(0):=0,
    \label{eq12}
\end{equation}
where $[\cdot]$ denotes the floor (entire) function. Thus, we define $a_1=a_1(x)=\left[\frac{1}{x}\right]$ and $a_n = a_1 (\tau ^{n-1}(x))$, $\forall n \in N$, where $\tau ^0(x)=x$, and $\tau^n(x)=\tau(\tau^{n-1}(x))$.
Then, from relation (\ref{eq12}), we have:
\[
    x= \frac{1}{\displaystyle a_1+\tau(x)} =
        \frac{1}{\displaystyle a_1+\frac{1}{\displaystyle a_2+\tau^2(x)}} = \ldots = [0;a_1,a_2,...,a_n+\tau^n(x)].
\]

The metrical theory of continued fractions expansions is about the sequence $(a_n)_{n\in N}$  of its incomplete quotients, and related sequences.
This theory started with Gauss' problem. In modern notation, Gauss wrote that
\begin{equation}
    \lim_{n\rightarrow\infty} \lambda\left(\left\{x\in[0,1); \tau^n(x)\leq z \right\}\right) = \frac{\log(z+1)}{\log 2}, \ 0\leq z\leq 1,
    \label{eq13}
\end{equation}
where $\lambda$ is the Lebesgue measure.  Gauss asked Laplace to
prove (\ref{eq13}) and to estimate the error-term $r_n(z)$,
defined by $r_n(z) := \lambda(\tau^{-n}([0,z])) -
\frac{\log(z+1)}{\log 2}$, $n\geq 1$. (Note that, when we omit the
logarithm base, we will consider the natural logarithm.) The first
one who proved (\ref{eq13}) and at the same time answered Gauss'
question was Kuzmin (1928), followed by L\' evy. From that time
on, a great number of such Gauss-Kuzmin theorems followed. To
mention a few: F. Schweiger (1968), P. Wirsing (1974), K.I.
Babenko (1978), and more recently by M. Iosifescu (1992).

Apart from regular continued fractions, there are many other
continued fractions expansions: Engel continued fractions, Rosen
expansions, the nearest integer continued fraction, the grotesque
continued fractions, etc.

\section{ANOTHER CONTINUED FRACTION EXPANSION}

We start this section by showing that any $x \in I_m :=[0, m-1]$, $m\in N_+$, $m\geq 2$, can be written in the form
\begin{equation}
    \frac{m^{-b_1(x)}}{\displaystyle 1+\frac{m^{-b_2(x)}}{1+\ddots}}:=[b_1(x), b_2(x), \ldots]_m
    \label{eq21}
\end{equation}
where $b_n = b_n(x)$ are integer values, belonging to the set $Z_{\geq-1}:=\{-1,0,1,2,\ldots\}$, for any $n\in N_+$.

{\bf Proposition 2.1} For any $x \in I_m:=[0,m-1]$, there exist integers numbers $b_n(x) \in \{-1,0,1,2,\ldots\}$ such that
\begin{equation}
    x = \frac{m^{-b_1(x)}}{\displaystyle 1+\frac{m^{-b_2(x)}}{1+\ddots}}
    \label{eq22}
\end{equation}
{\bf Proof}. If $x\in [0, m-1]$, then we can find an integer $b_1(x)\in Z_{\geq-1}$ such that
\begin{equation}
    \frac{1}{m^{b_1(x)+1}}<x<\frac{1}{m^{b_1(x)}}.
    \label{eq23}
\end{equation}
Thus, there is a $p\in[0,1]$ such that
\[
    x = (1-p)\frac{1}{m^{b_1(x)}} + p\frac{1}{m^{b_1(x)+1}} = \frac{m-mp+p}{m} m^{-b_1(x)}.
\]
If we set $x_1 = \frac{mp-m}{m-mp+p}$, then we can write $x$ as
\[
    x = \frac{m^{-b_1(x)}}{1+x_1}.
\]
Since $x_1\in[0,m-1]$, we can repeat the same iteration and obtain
\[
    x = \frac{m^{-b_1(x)}}{\displaystyle 1+\frac{m^{-b_2(x)}}{1+\ddots}}
\]
which completes the proof.

Next, we define on $I_m:=[0,m-1]$, $m\in N_+$, $m\geq 2$, the transformation $\tau_m$ by
\[
    \tau_m: I_m \rightarrow I_m,
\]
\begin{equation}
    \tau_m(x):=m^{\frac{\log x^{-1}}{\log m}-\left[\frac{\log x^{-1}}{\log m}\right]}-1, x\neq 0, \tau(0):=0,
    \label{eq24}
\end{equation}
where $[\cdot]$ denotes the floor (entire) function.

For any $x \in I_m$, put
\[
    b_n(x)=b_1\left(\tau^{n-1}_m(x)\right), n \in N_+,
\]
\[
    b_1(x) = \left[\frac{\log x^{-1}}{\log m}\right], x \neq 0, b_1(0)=\infty.
\]
Let $\Omega_m$ be the set of all irrational numbers from $I_m$. In
the case when $x\in I_m\backslash \Omega_m$, we have
\[
    b_n(x) = \infty, \forall n\geq k(x) \geq m,\mbox{ and } b_n(x) \in Z_{\geq -1}, \forall n < k(x).
\]
Therefore, in the rational case, the continued fraction expansion (\ref{eq21}) is finite, unlike the irrational case,
when we have an infinite number of incomplete quotients from the set $\{-1,0,1,2,\ldots\}$.

Let $\omega \in \Omega_m$. We have
\[
    \omega = m^{\log_m \omega} = m ^{-\frac{\log \omega^{-1}}{\log m}} = \frac{m^{-\left[\frac{\log \omega^{-1}}{\log m}\right]}}{m^{\frac{\log \omega^{-1}}{\log m}-\left[\frac{\log \omega^{-1}}{\log m}\right]}} = \frac {m^{-b_1(\omega)}}{1+\tau_m(\omega)}.
\]
Since,
\begin{eqnarray}
    \tau_m(\omega) &=& m^{\log_m\tau_m(\omega)} = m^{-\frac{\log \tau_m^{-1}(\omega)}{\log m}} = \frac{m^{-\left[\frac{\log \tau_m^{-1}(\omega)}{\log m}\right]}}{m^{\frac{\log \tau_m^{-1}(\omega)}{\log m} - \left[\frac{\log \tau_m^{-1}(\omega)}{\log m}\right]}} \nonumber \\
&=&\frac{m^{-b_1(\tau_m(\omega))}}{1+\tau_m(\tau_m(\omega))} = \frac{m^{-b_2(\omega)}}{1+\tau^2_m(\omega)} \nonumber
\end{eqnarray}
then, we have
\begin{equation}
    \omega = \frac{m^{-b_1(\omega)}}{\displaystyle 1+ \frac{m^{-b_2(\omega)}}{1+\tau^2_m(\omega)}} = \ldots = \frac{m^{-b_1(\omega)}}{\displaystyle 1+ \frac{m^{-b_2(\omega)}}{\displaystyle 1+ \ddots + \frac{m^{-b_n(\omega)}}{1+ \tau^n_m(\omega)}}}
    \label{eq25}
\end{equation}
If $[b_1(\omega)] = m^{-b_1(\omega)}$, and $[b_1(\omega), b_2(\omega), \ldots, b_n(\omega)] = \frac{m^{-b_1(\omega)}}{1+[b_2(\omega),b_3(\omega),\ldots,b_n(\omega)]}$, $\forall n\geq 2$, then (\ref{eq25}) can be written as
\begin{eqnarray}
    \omega &=& \left[b_1(\omega) + \frac{\log(1+\tau_m(\omega))}{\log m}\right] = \left[b_1(\omega), b_2(\omega) + \frac{\log(1+\tau^2_m(\omega))}{\log m}\right] = \ldots = \nonumber\\
    &=& \left[b_1(\omega), b_2(\omega), \ldots, b_{n-1}(\omega), b_n(\omega) + \frac{\log(1+\tau^n_m(\omega))}{\log m}\right].\nonumber
\end{eqnarray}

It is obvious that we have the relations
\begin{equation}
    \tau_m(\omega) = \frac{m^{-b_2(\omega)}}{1+\tau^2_m(\omega)}, \ldots, \tau^{n-1}_m(\omega) = \frac{m^{-b_n(\omega)}}{1+\tau^n_m(\omega)}, \forall n\in N_+, \forall \omega \in \Omega_m,
    \label{eq26}
\end{equation}

\section{CONVERGENTS. BASIC PROPERTIES}
In this section we define and give the basic properties of the
convergents of this continued fraction expansion.

{\bf Definition 3.1} A finite truncation in (\ref{eq21}), i.e.
\begin{equation}
    \omega_n(\omega):=\frac{p_n(\omega)}{q_n(\omega)} = [b_1(\omega), b_2(\omega), \ldots, b_n(\omega)]_m, n\in N_+
    \label{eq31}
\end{equation}
is called the $n$-th convergent of $\omega$.

The integer valued functions sequences $(p_n)_{n\in N}$ and $(q_n)_{n\in N}$ can be recursively defined by the formulae:
\begin{eqnarray}
    p_n(\omega) &=& m^{b_n(\omega)}p_{n-1}(\omega) + m^{b_{n-1}(\omega)}p_{n-2}, \forall n\geq 2, \nonumber \\
    q_n(\omega) &=& m^{b_n(\omega)}q_{n-1}(\omega) + m^{b_{n-1}(\omega)}q_{n-2}, \forall n\geq 2,
    \label{eq32}
\end{eqnarray}
with $p_0(\omega)=0$, $q_0(\omega) = 1$, $p_1(\omega)=1$ and $q_1(\omega)=m^{b_1(\omega)}$.

By induction, it is easy to prove that
\begin{equation}
    p_n(\omega)q_{n+1}(\omega) - p_{n+1}(\omega)q_n(\omega) = (-1)^{n+1}m^{b_1(\omega) + \ldots + b_n(\omega)}, \forall n \in N_+,
    \label{eq33}
\end{equation}
and that
\begin{equation}
    \frac{m^{-b_1(\omega)}}{\displaystyle 1+ \frac{m^{-b_2(\omega)}}{\displaystyle 1+ \ddots + \frac{m^{-b_n(\omega)}}{1+ t}}} = \frac{p_n(\omega) + tm^{b_n(\omega)}p_{n-1}(\omega)}{q_n(\omega) + tm^{b_n(\omega)}q_{n-1}(\omega)}, \forall n \in N_+, t\geq 0.
    \label{eq34}
\end{equation}

Now, combining the relations (\ref{eq25}) and (\ref{eq34}), it can be shown that, for any $\omega \in \Omega_m$, we have
\begin{equation}
    \omega = \frac{p_n(\omega) + \tau_m^n(\omega)m^{b_n(\omega)}p_{n-1}(\omega)}{q_n(\omega) + \tau_m^n(\omega)m^{b_n(\omega)}q_{n-1}(\omega)}, \forall n \in N_+.
    \label{eq35}
\end{equation}

\section{MAIN RESULT}
At this moment, we are able to present the main result of the paper, which is the convergence of this new continued fraction expansion, i.e. we must show that
\[
    \omega = \lim_{n\rightarrow\infty}[b_1(\omega), b_2(\omega), \ldots, b_n(\omega)]_m,
\]
for any $\omega \in \Omega_m$.

{\bf Theorem} For any $\omega \in \Omega_m:=I_m \backslash Q$, we have
\begin{equation}
    \omega - \omega_n(\omega) = \frac{(-1)^n\tau^n_m(\omega)m^{b_1(\omega)+\ldots+ b_n(\omega)}}{q_n(\omega)\left(q_n(\omega)+\tau^n_m(\omega)m^{b_n(\omega)}q_{n-1}(\omega)\right)}, \forall n \in N_+.
    \label{eq41}
\end{equation}
For any $\omega \in \Omega_m$, we have
\[
    \frac{m^{b_1(\omega)+\ldots+b_n}(\omega)}{q_n(\omega)\left(q_{n+1}(\omega)+(m-1)^{n+1}m^{b_{n+1}(\omega)}q_{n}(\omega)\right)} <
    |\omega - \omega_n(\omega)|<
\]
\begin{equation}
    < \frac{1}{\max(F_n, m^{b_1(\omega)+\ldots+b_n(\omega)})}, \forall n\in N_+,
    \label{eq42}
\end{equation}
and
\begin{equation}
    \lim_{n\rightarrow\infty}\omega_n(\omega) = \omega
    \label{eq43}
\end{equation}
Here $F_n$ denotes the $n$-th Fibonacci number.\\
{\bf Proof.} Using relations (\ref{eq33}) and (\ref{eq35}), we obtain:
\begin{eqnarray}
    \omega - \omega_n(\omega) &=& \frac{p_n(\omega) + \tau^n_m(\omega) m^{b_n(\omega)}p_{n-1}(\omega)}{q_n(\omega) + \tau^n_m(\omega)m^{b_n(\omega)}q_{n-1}(\omega)} - \frac{p_n(\omega)}{q_n(\omega)} \nonumber \\
    &=& \frac{(-1)^n\tau^n_m(\omega)m^{b_1(\omega)+\ldots+ b_n(\omega)}}{q_n(\omega)\left(q_n(\omega)+\tau^n_m(\omega)m^{b_n(\omega)}q_{n-1}(\omega)\right)}. \nonumber
\end{eqnarray}

Next, by (\ref{eq26}) and (\ref{eq41}), it follows:
\begin{eqnarray}
    |\omega - \omega_n(\omega)| & =& \frac{\tau^n_m(\omega)m^{b_1(\omega)+\ldots+ b_n(\omega)}}{q_n(\omega)\left(q_n(\omega)+ \tau^n_m(\omega)m^{b_n(\omega)}q_{n-1}(\omega)\right)} \nonumber \\
    &=& \frac{m^{-b_{n+1}(\omega)}}{1+\tau^{n+1}_m(\omega)} \cdot \frac{m^{b_1(\omega)+\ldots+ b_n(\omega)}} {q_n(\omega)\left(q_n(\omega)+ \frac{m^{-b_{n+1}(\omega)}}{1+\tau^{n+1}_m(\omega)} m^{b_n(\omega)}q_{n-1}(\omega)\right)}   \nonumber \\
    &=& \frac{m^{b_1(\omega)+\ldots+ b_n(\omega)}} {m^{b_{n+1}(\omega)}q_n(\omega)\left(q_n(\omega) + \tau^{n+1}_m(\omega)q_n(\omega) + m^{-b_{n+1}(\omega)}m^{b_n(\omega)}q_{n+1}(\omega)\right)} \nonumber \\
    &=& \frac{m^{b_1(\omega)+\ldots+ b_n(\omega)}} {q_n(\omega)\left(q_{n+1}(\omega) + \tau^{n+1}_m(\omega)m^{b_{n+1}(\omega)}q_n(\omega)\right)}
    \label{eq44}
\end{eqnarray}
Now, we know that the Fibonacci numbers are defined by recurrence
\[
    F_{n+1} = F_{n} + F_{n-1}, \forall n \in N_+, \mbox{ and } F_0 = F_1=1.
\]
Also, from the recurrence relation (\ref{eq32}), we infer that
\begin{equation}
    p_{n+1} \geq F_{n+1} \mbox{ and } q_n \geq F_n, \forall n \in N_+, n \geq 2.
    \label{eq45}
\end{equation}
Also, we have that
\begin{eqnarray}
    q_n(\omega) &=& m^{b_n(\omega)} q_{n-1}(\omega) + m^{b_{n-1}(\omega)}q_{n-2}(\omega) \geq m^{b_n(\omega)}q_{n-1}(\omega) \geq \nonumber\\
    & \geq& m^{b_n(\omega)} m^{b_{n-1}(\omega)} q_{n-2}(\omega) \geq \ldots \geq m^{b_1(\omega)+\ldots+b_n(\omega)}q_0(\omega). \nonumber
\end{eqnarray}
i.e.
\begin{equation}
    q_n(\omega) \geq m^{b_1(\omega)+\ldots+b_n(\omega)}, \forall n \in N_+.
    \label{eq46}
\end{equation}
Thus, from relations (\ref{eq45}) and (\ref{eq46}), we have that
\[
    q_n(\omega) \geq \max(F_n, m^{b_1(\omega)+\ldots+b_n(\omega)}), \forall n \in N_+.
\]
Now, since the transformation $\tau_m$ belonging to $(0, m-1)$ and
from the last two relations, we can show that
\[
    \frac{m^{b_1(\omega)+\ldots+ b_n(\omega)}}{q_n(\omega)\left(q_{n+1}(\omega) + \tau^{n+1}_m(\omega) m^{b_{n+1}(\omega)}
    q_n(\omega)\right)} \leq \frac{m^{b_1(\omega)+\ldots+ b_n(\omega)}}{q_n(\omega)q_{n+1}(\omega)}\leq
\]
\[
 \leq \frac{1}{q_n(\omega)} \leq \frac{1}{\max(F_n, m^{b_1(\omega)+\ldots+b_n(\omega)})}.
\]
It is obvious that the left inequality is true. Since $\max(F_n,
m^{b_1(\omega)+\ldots+b_n(\omega)})$ is an increasing function, we
have
\[
    \lim_{n \rightarrow \infty}\omega_n(\omega) = \omega.
\]
The proof is complete.

\section{REMARK}
This paper is the first one which addresses this type of continued
fraction expansion, and will be followed by other papers which
will present the metrical theory of this expansion, the principal
aim being solving Gauss' problem.

\end{document}